\documentclass[letterpaper,11pt]{article}
\usepackage{amsmath, amsthm, amssymb}    
\usepackage{bridges}			
\usepackage{graphicx}			
\usepackage[colorlinks=true, urlcolor=blue, citecolor=black, linkcolor=black]{hyperref}  
\usepackage{subfloat}
\usepackage{subcaption}			

\usepackage[utf8]{inputenc}
\usepackage{amsfonts}
\usepackage{mathtools}
\usepackage{wrapfig}
\usepackage{tikz}

\title{Surfaces in The Tesseract}

\author{Manuel Est\'evez\textsuperscript{1}, \'Erika Rold\'an\textsuperscript{2,3}, and Henry Segerman\textsuperscript{4}
\vspace{10pt}\\
\textsuperscript{1} ScaDS.AI, Leipzig University; estevez@mis.mpg.de\\
\textsuperscript{2} ScaDS.AI, Leipzig University; roldan@mis.mpg.de\\
\textsuperscript{3} Max Planck Institute for Mathematics in the Sciences; roldan@mis.mpg.de\\
\textsuperscript{4} Oklahoma State University; henry@segerman.org
} 

\date{}					

\begin{document}

\maketitle

\thispagestyle{empty}

\begin{abstract}
How can we visualize all the surfaces that can be made from the faces of the tesseract? In recent work, Aveni, Govc, and Rold\'an showed that the torus and the sphere are the only closed surfaces that can be realized by a subset of two-dimensional faces of the tesseract. They also gave an exhaustive list of all the isomorphic types of embedings of these two surfaces. Here, we generate 3D models of all these surfaces. We also exhibit, with the help of some hyper-ants, the minimum realization of the M\"obius strip on the tesseract. 
\end{abstract}

\section{The hypercube and its group of symmetries}
\label{background}
We denote the $n$-dimensional unit cube by $Q^n=[0,1]^n$, and its set of vertices by $Q_0^n$. Each vertex of $Q^n$ can be represented by an element of the set of all $n$-tuples with binary entries $\{0,1\}^n$. We denote by $Q_1^n$ the one-dimensional skeleton of $Q^n$, that is, the set of its vertices and edges. We observe that $Q_1^n$ is the graph with vertex set $Q_0^n$ and an edge between two vertices if and only if they differ in exactly one coordinate. Similarly, $Q_2^n$ denotes the two-dimensional skeleton of $Q^n$, which consists of the one-dimensional skeleton plus all its two-dimensional faces (we will refer to these as simply \emph{faces}).
We refer to a subset of $Q^n_2$ as a \emph{two-dimensional cubical complex}.

\begin{wrapfigure}[20]{l}{0.6\textwidth}
\centering
\vspace{-9pt}
\tikzset{unode/.style = {
    circle, 
    draw=cyan!20!black, 
    thick,
    fill={rgb:black,1;white,4},
    inner sep=2.3pt,
    minimum size=2.3pt } }
\tikzset{uedge/.style = {
    draw=cyan!20!black, 
    very thick} }

\begin{tikzpicture} [thick,scale=.95, every node/.style={scale=0.6}]
    \foreach \i in {0,1,2,3,...,7}{ 
        \path[uedge] (\i*45:4) -- (45+\i*45:4) {};
        \path[uedge] (\i*45:1.646) -- (45+\i*45:4){};
        \path[uedge] (\i*45:1.646) -- (-45+\i*45:4){};
        \path[uedge] (\i*45:1.646) -- (135+\i*45:1.646){};
        }
        {
        \path[uedge] (0*45:4) -- (45+0*45:4)  node[pos=.5, rotate=-70, above, scale=1.5, align=right]{(*,0,1,0)};
        }
        \node[unode] at (0*45:1.646) {(0,0,0,0)};
        \node[unode] at (5*45:1.646) {(0,1,0,0)};
        \node[unode] at (6*45:4) {(1,1,0,0)};
        \node[unode] at (7*45:4) {(1,0,0,0)};

        \node[unode] at (1*45:4) {(0,0,1,0)};
        \node[unode] at (2*45:1.646) {(0,1,1,0)};
        \node[unode] at (7*45:1.646) {(1,1,1,0)};
        \node[unode] at (0*45:4) {(1,0,1,0)};
        
        \node[unode] at (3*45:1.646) {(0,0,0,1)};
        \node[unode] at (4*45:4) {(0,1,0,1)};
        \node[unode] at (5*45:4) {(1,1,0,1)};
        \node[unode] at (6*45:1.646) {(1,0,0,1)};

        \node[unode] at (2*45:4) {(0,0,1,1)};
        \node[unode] at (3*45:4) {(0,1,1,1)};
        \node[unode] at (4*45:1.646) {(1,1,1,1)};
        \node[unode] at (1*45:1.646) {(1,0,1,1)};
    
\end{tikzpicture}
\caption{Vertices and edges of $Q^4$.}
\label{hcube}
\end{wrapfigure}
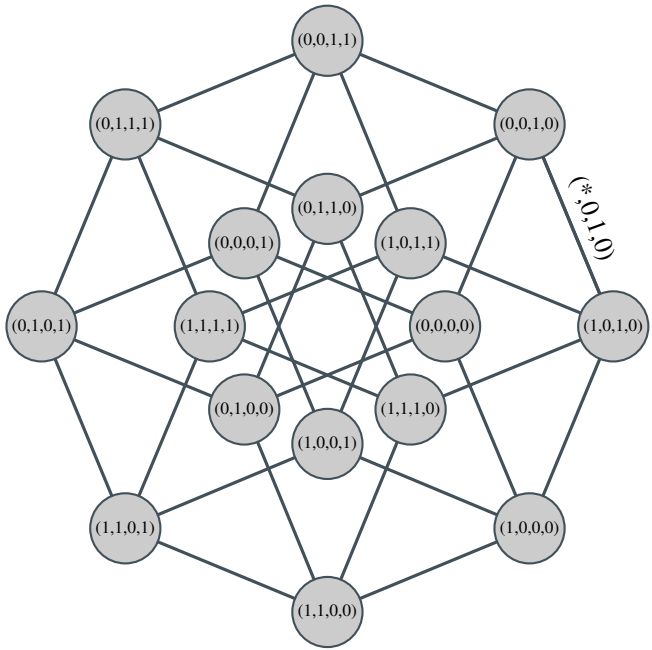
For a unit cube $Q^3=[0,1]^3$, consider its sets of vertices $Q^3_0$, edges $Q^3_1$, faces $Q^3_2$, and of course the set consisting of the cube itself $Q^3_3$. The elements of these sets are the \textit{cells} of $Q^3$. Geometrically, every cell of $Q^n$ is a product of vertices and intervals, and therefore can be encoded combinatorially as an element of $\{0,1,*\}^n$. Here a * in an entry implies that in the product, the whole interval $I$ is considered in that direction. For example, in the tesseract, $Q^4$, the four-tuple $(*,0,1,0)$ encodes the edge displayed in the upper right part of Figure \ref{hcube}.
Thus, every subcomplex of $Q^n$ can be represented as a subset of $\{0,1,*\}^n$. We call this \textit{star notation}. 
We can use this notation to easily calculate the number of cells in $Q^4$. This is the size of $\{0,1,*\}^4$, yielding $3^4=81$ cells. 

The $n$-dimensional \textit{hyperoctahedral group}, denoted by $B_n$ is a group of order $2^nn!$ describing the symmetries of an $n$-dimensional hypercube. The group $B_2$ is precisely the group of symmetries of the square: the dihedral group of order eight. However, $B_n$ can also be regarded as the \textit{group of signed permutation matrices}, which consists of all $n \times n$ invertible matrices with entries $\{1,-1,0\}$, with matrix multiplication as its group operation. Since we are interested in the symmetries of the tesseract, $B_4$ will play an important role when classifying all of its surfaces.

\section{Cubical surfaces}
\label{CubicalSurfaces}

\begin{figure}[htb]
\vspace{8pt}
    \centering
\subfloat[A]{
        \includegraphics[height=1.02in]{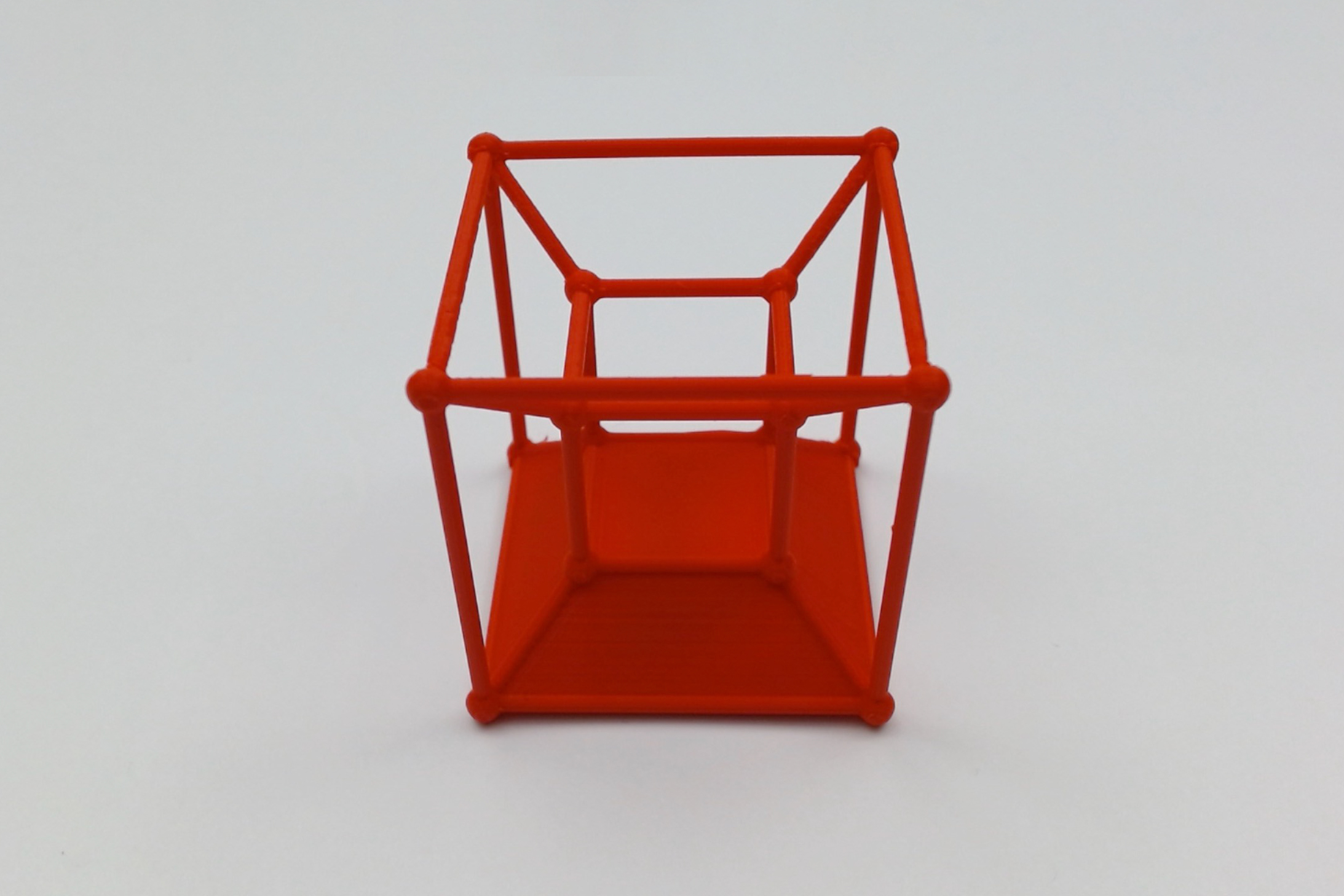}
    }
\subfloat[A]{
        \includegraphics[height=1.02in]{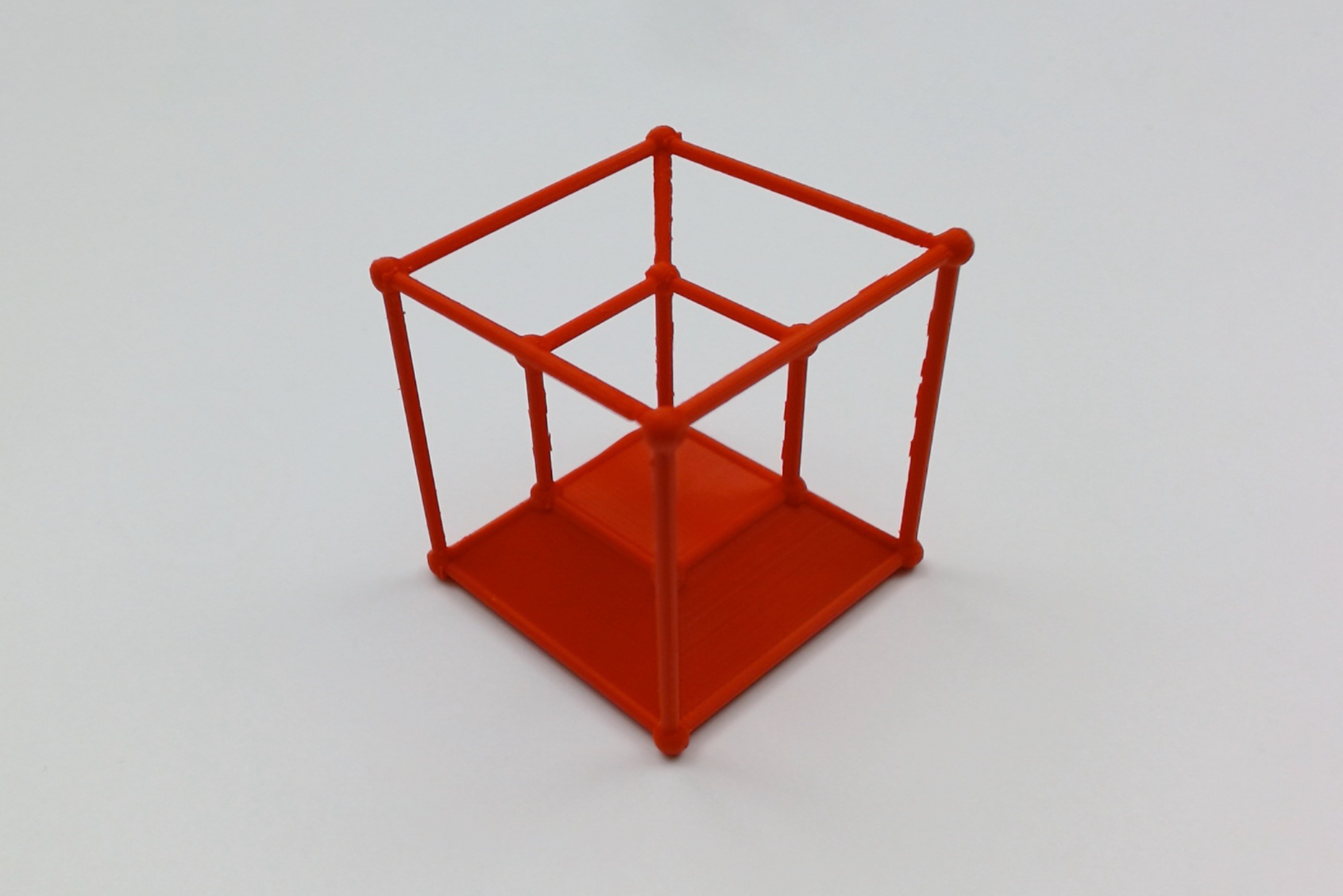}
    }
\subfloat[B]{
        \includegraphics[height=1.02in]{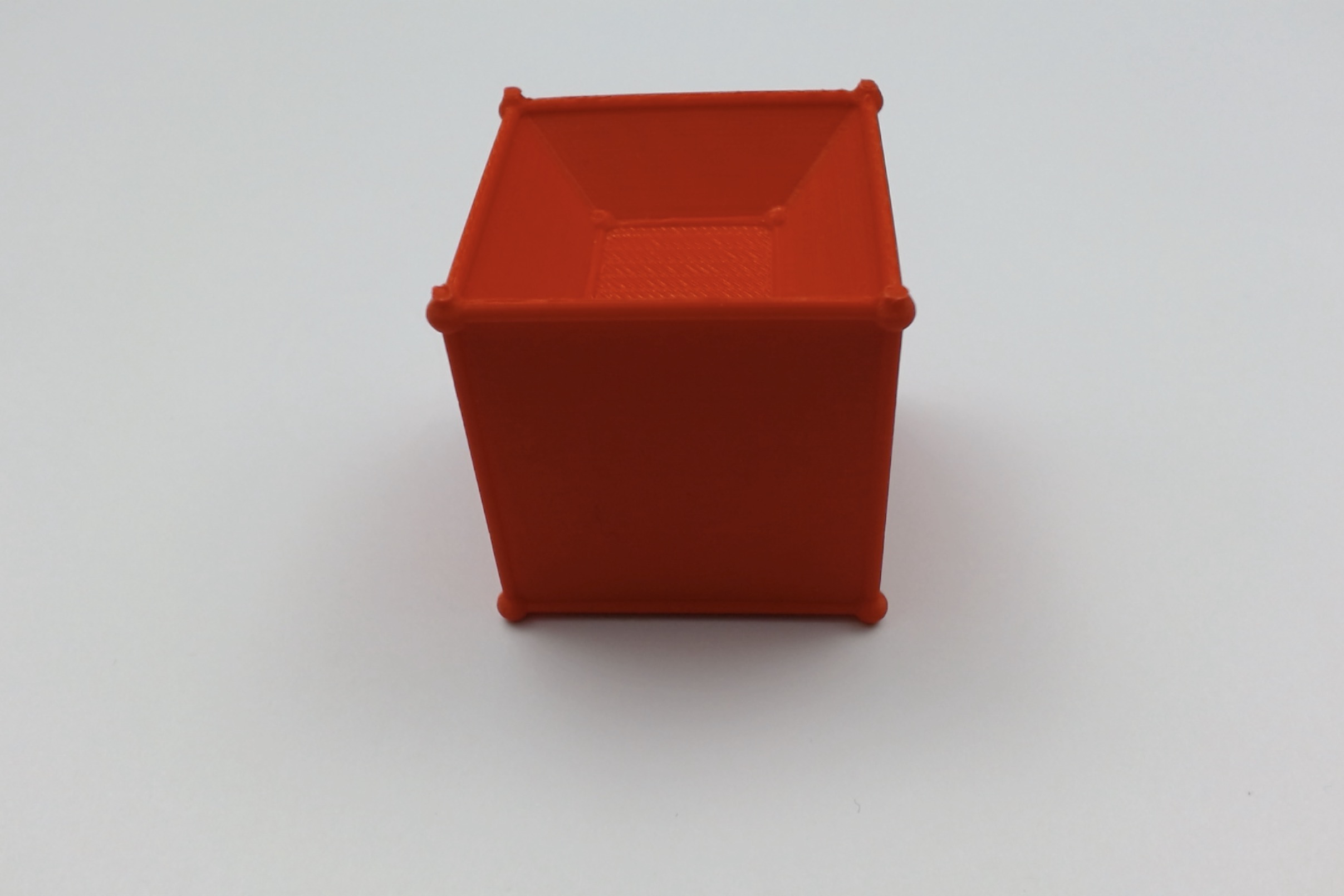}
    }
\subfloat[B]{
        \includegraphics[height=1.02in]{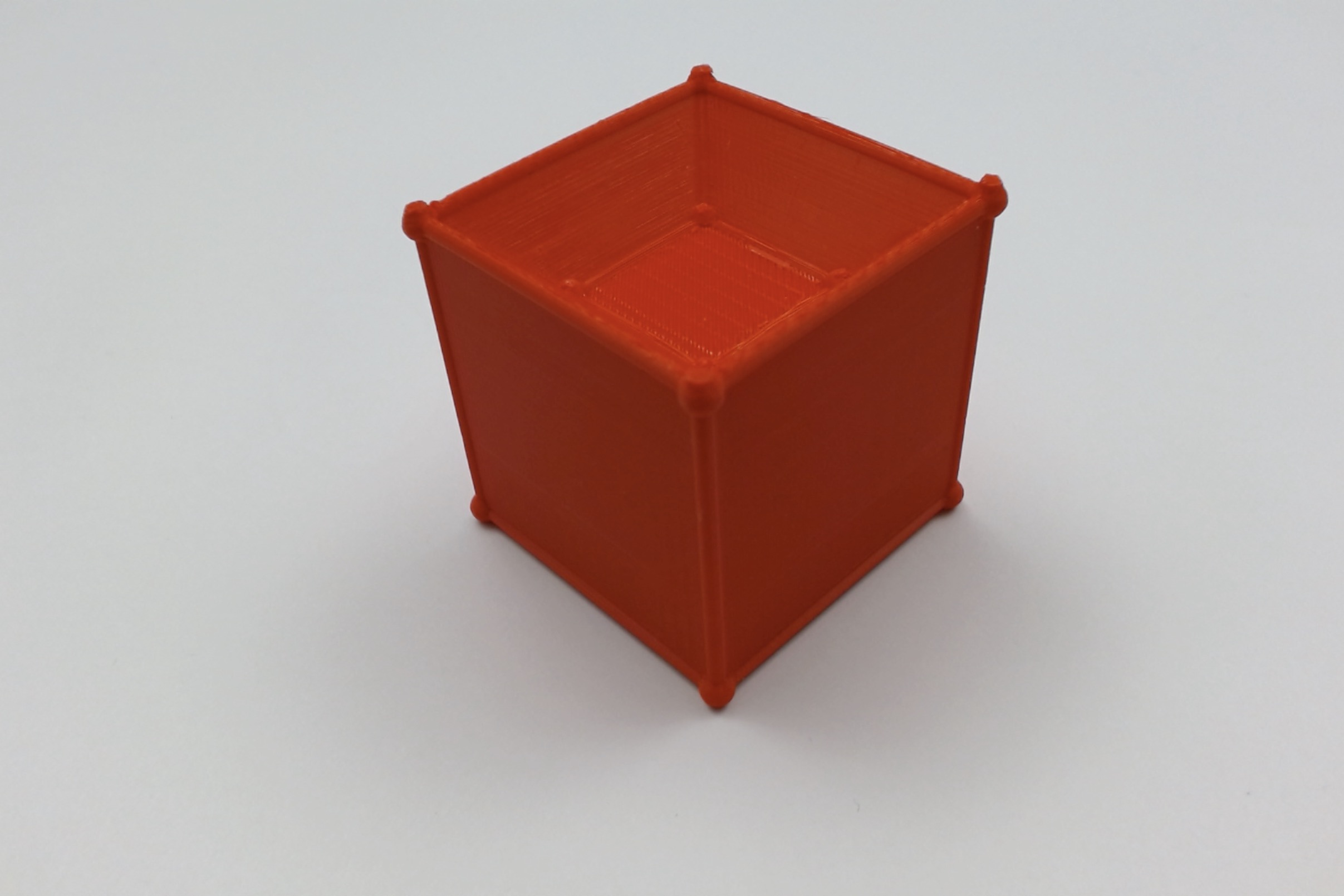}
    }

\subfloat[C]{
        \includegraphics[height=1.02in]{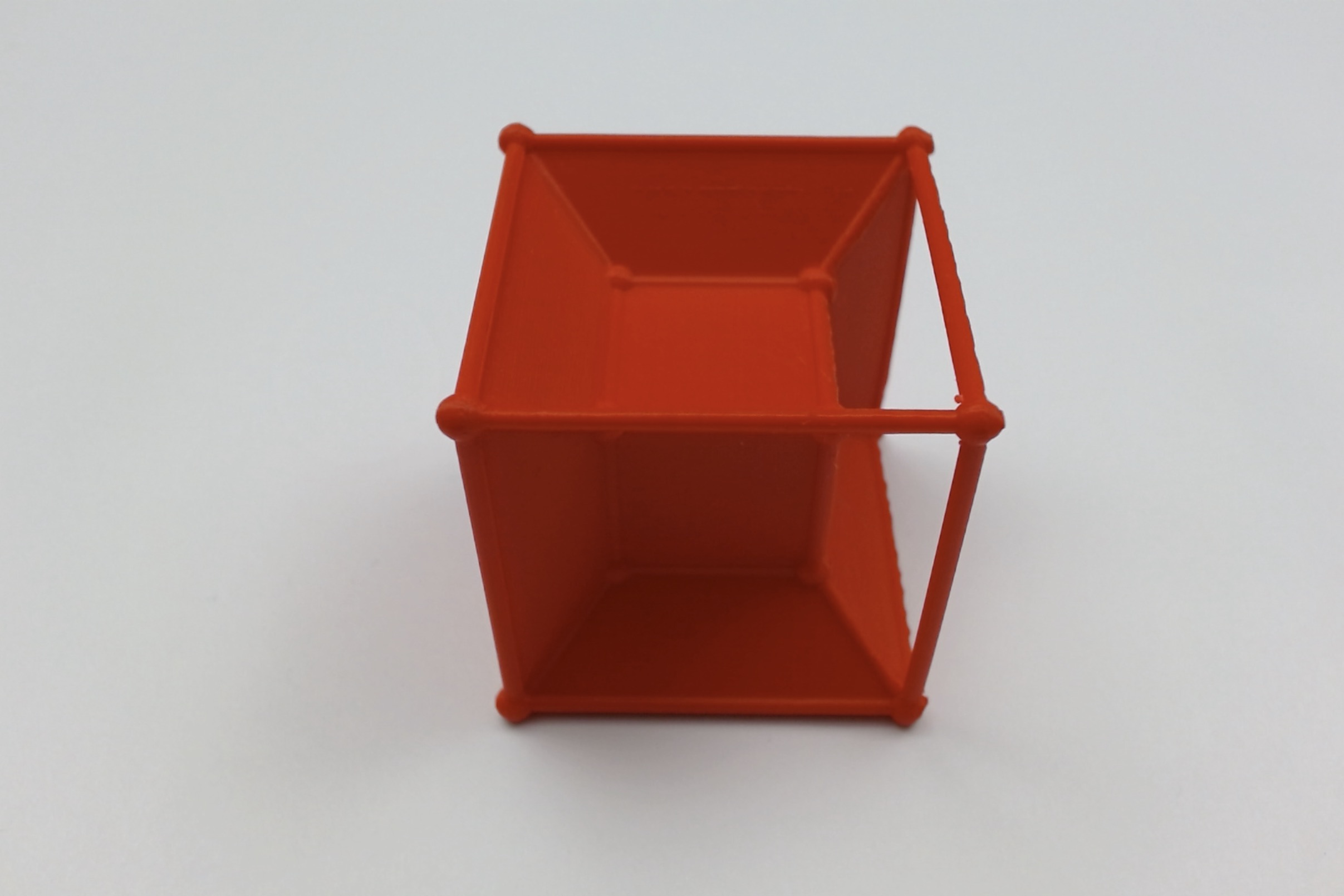}
    }
\subfloat[C]{
        \includegraphics[height=1.02in]{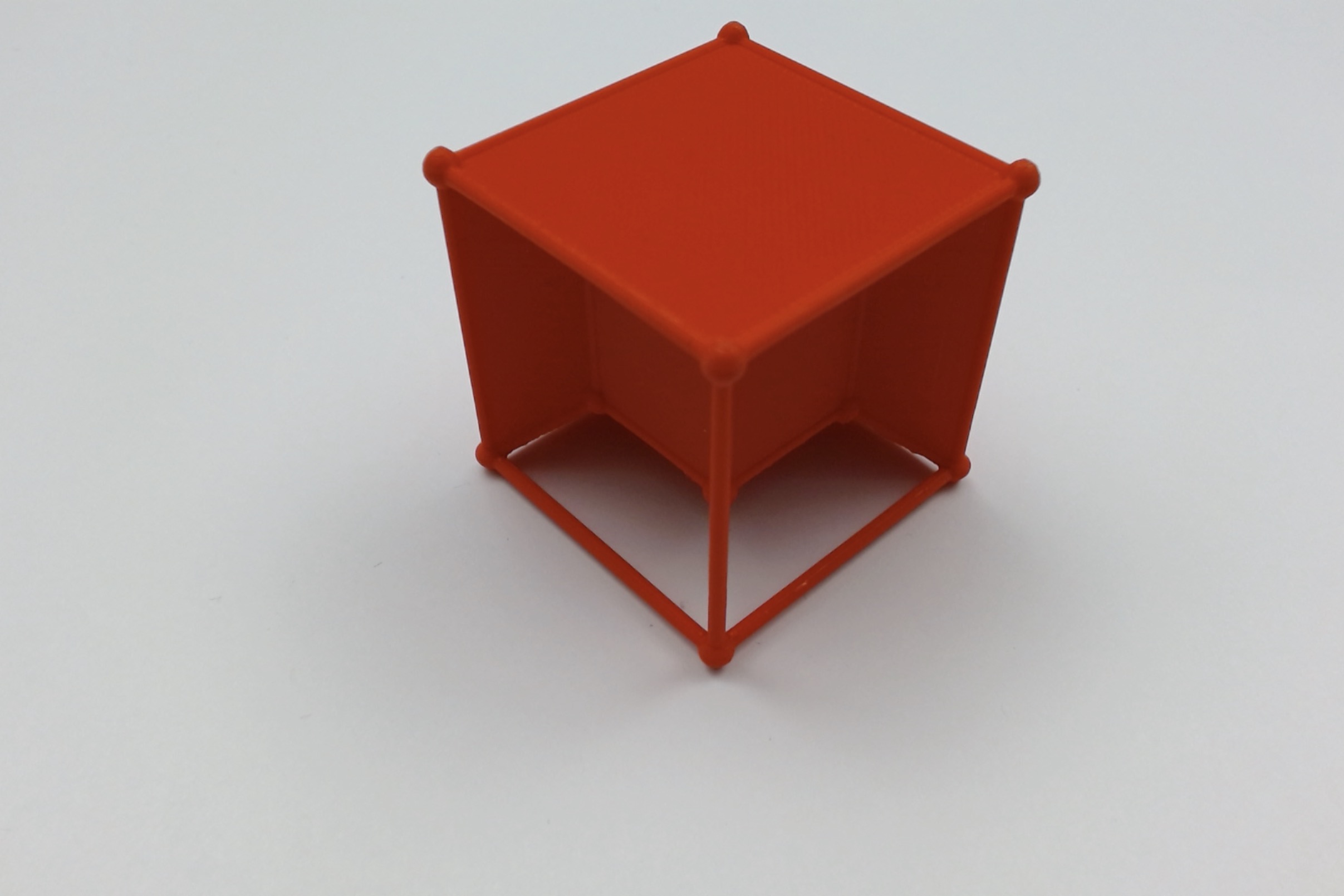}
    }
\subfloat[D]{
        \includegraphics[height=1.02in]{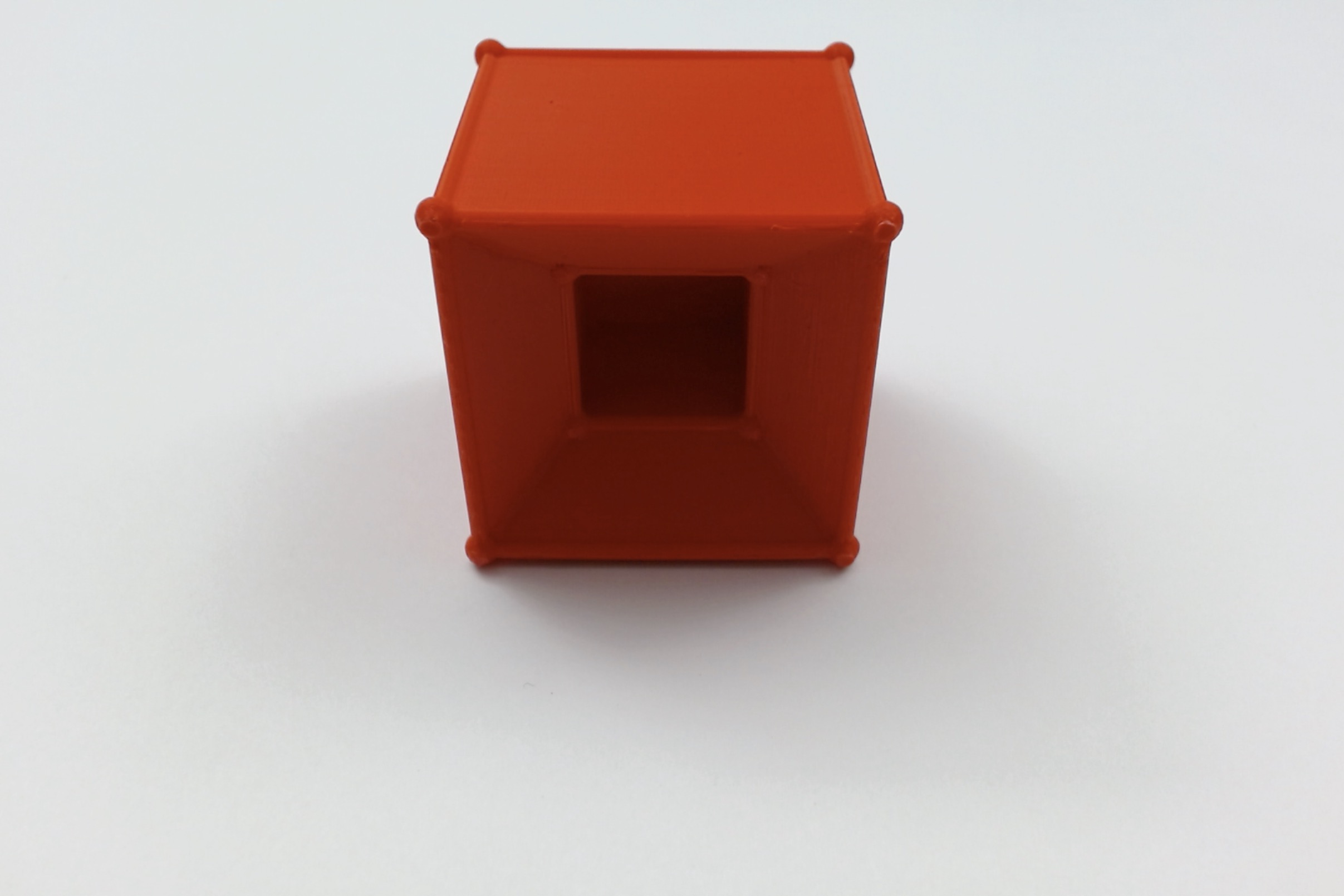}
    }
\subfloat[D]{
        \includegraphics[height=1.02in]{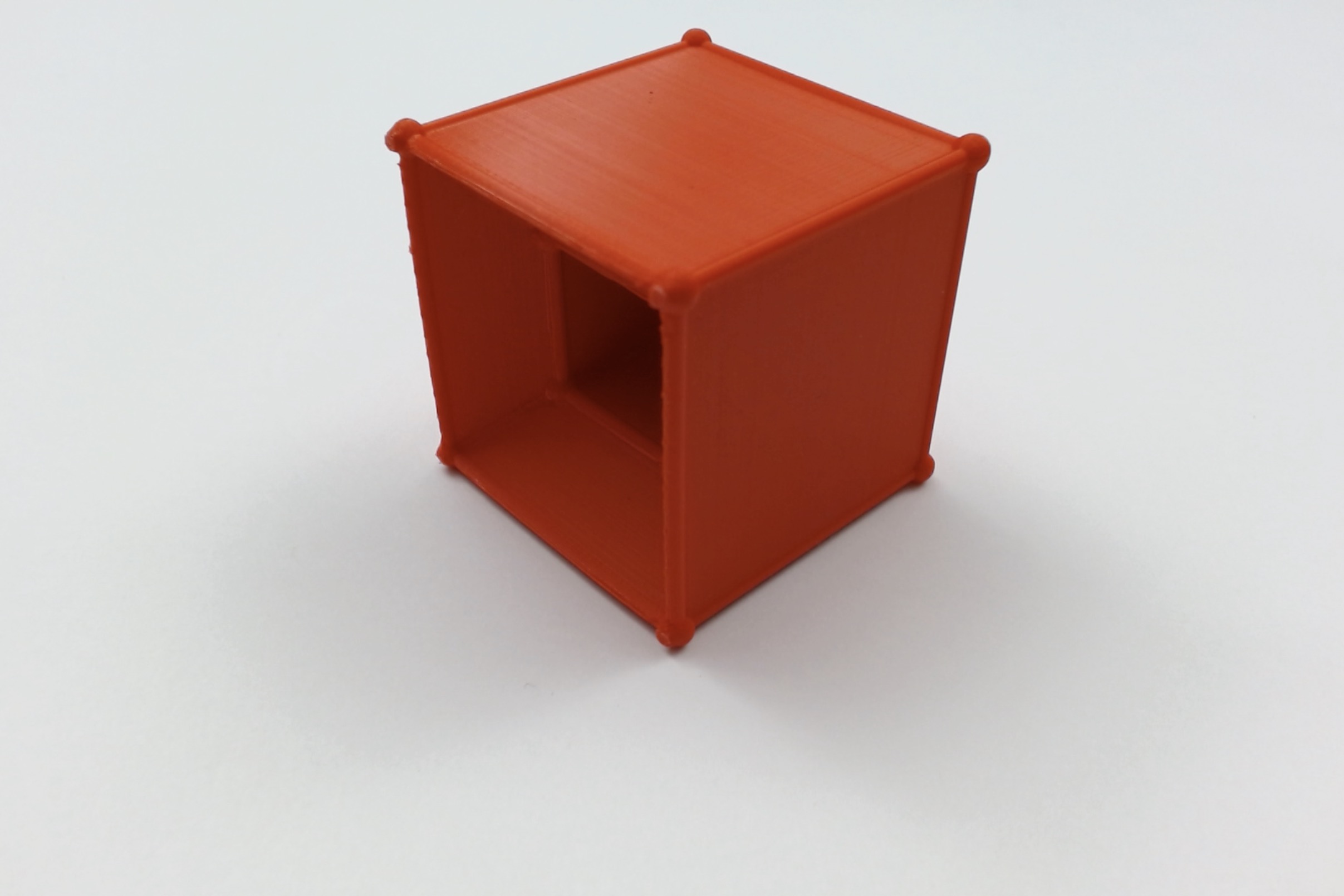}
    }

\subfloat[E]{
        \includegraphics[height=1.02in]{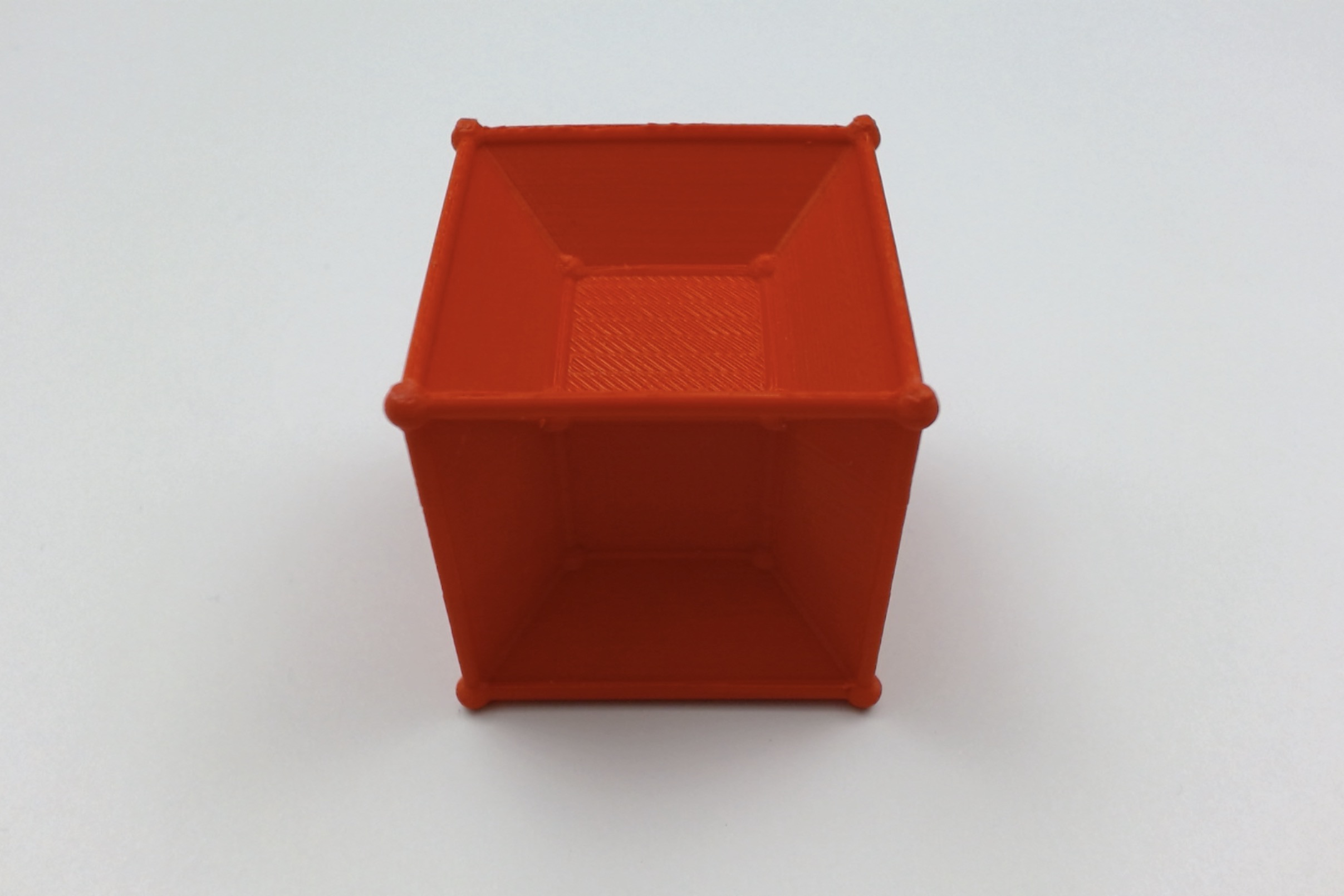}
    }
\subfloat[E]{
        \includegraphics[height=1.02in]{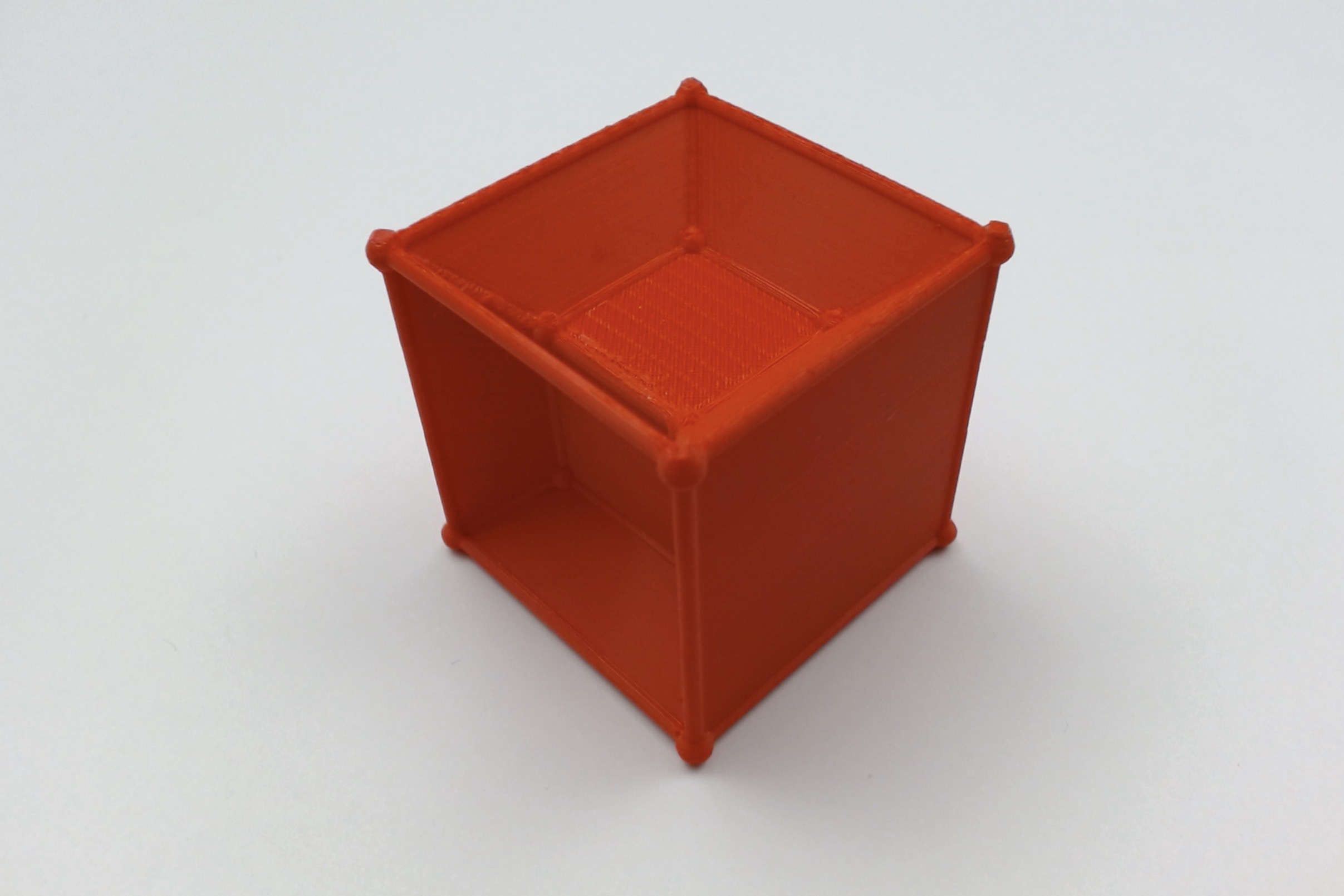}
    }
\subfloat[F]{
        \includegraphics[height=1.02in]{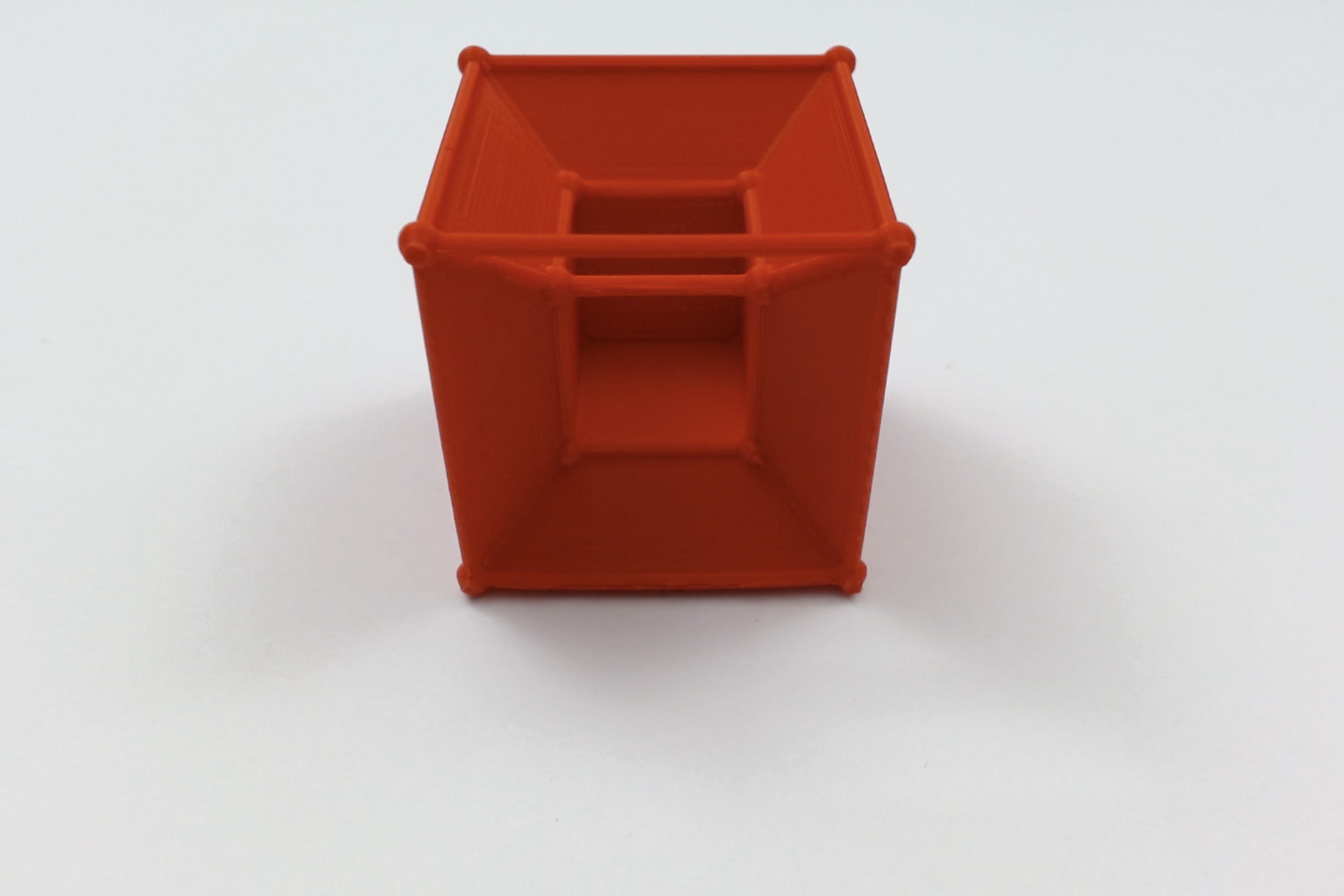}
    }
\subfloat[F]{
        \includegraphics[height=1.02in]{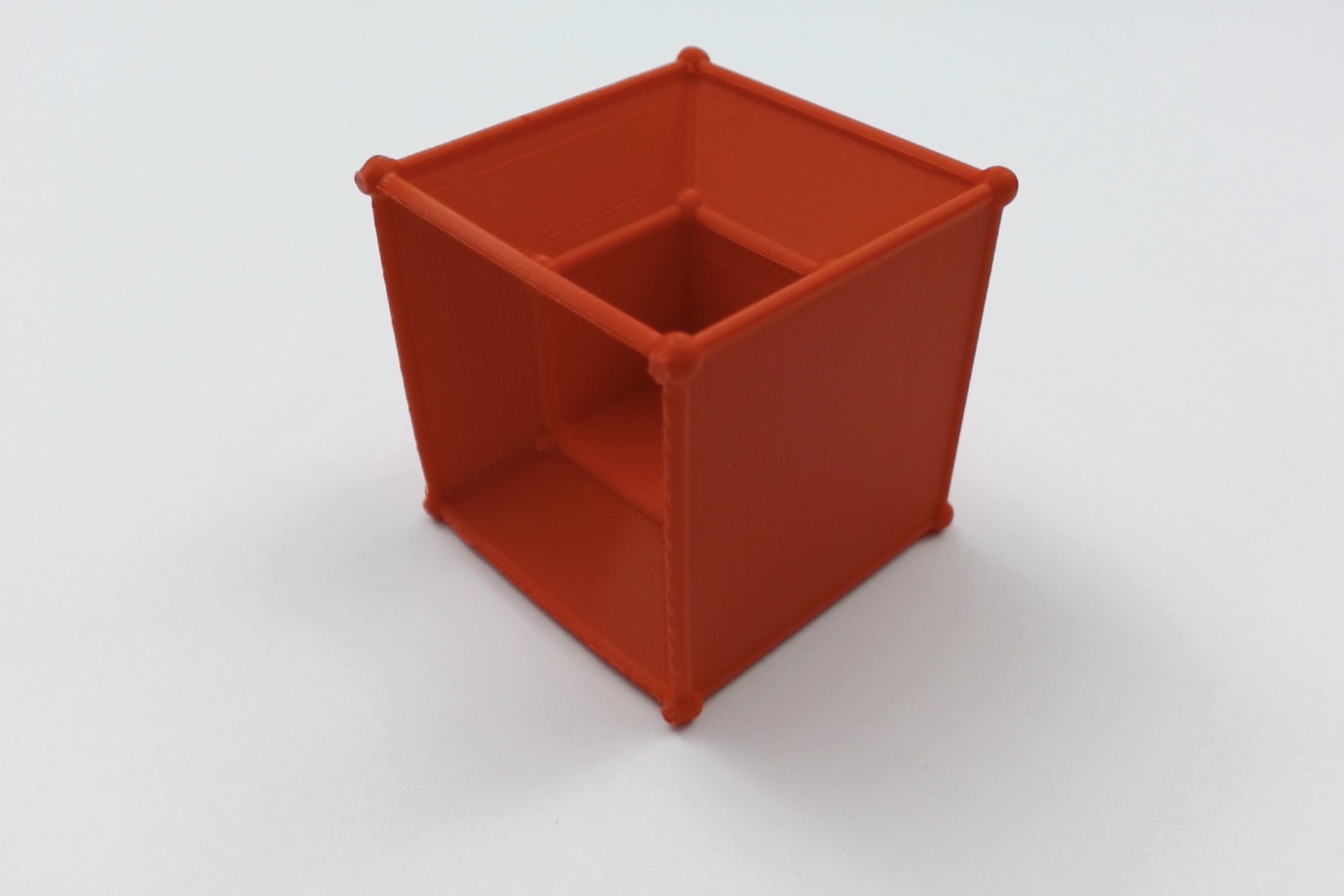}
    }
\caption{Embeddings of $\mathbb{S}^2$ in the tesseract. The surfaces $A$ through $F$ have 8, 24, 8, 24, 32, and 24 possible representatives respectively.}
\label{prints}
\end{figure}

\begin{figure}
    \centering
        \includegraphics[height=1.2in]{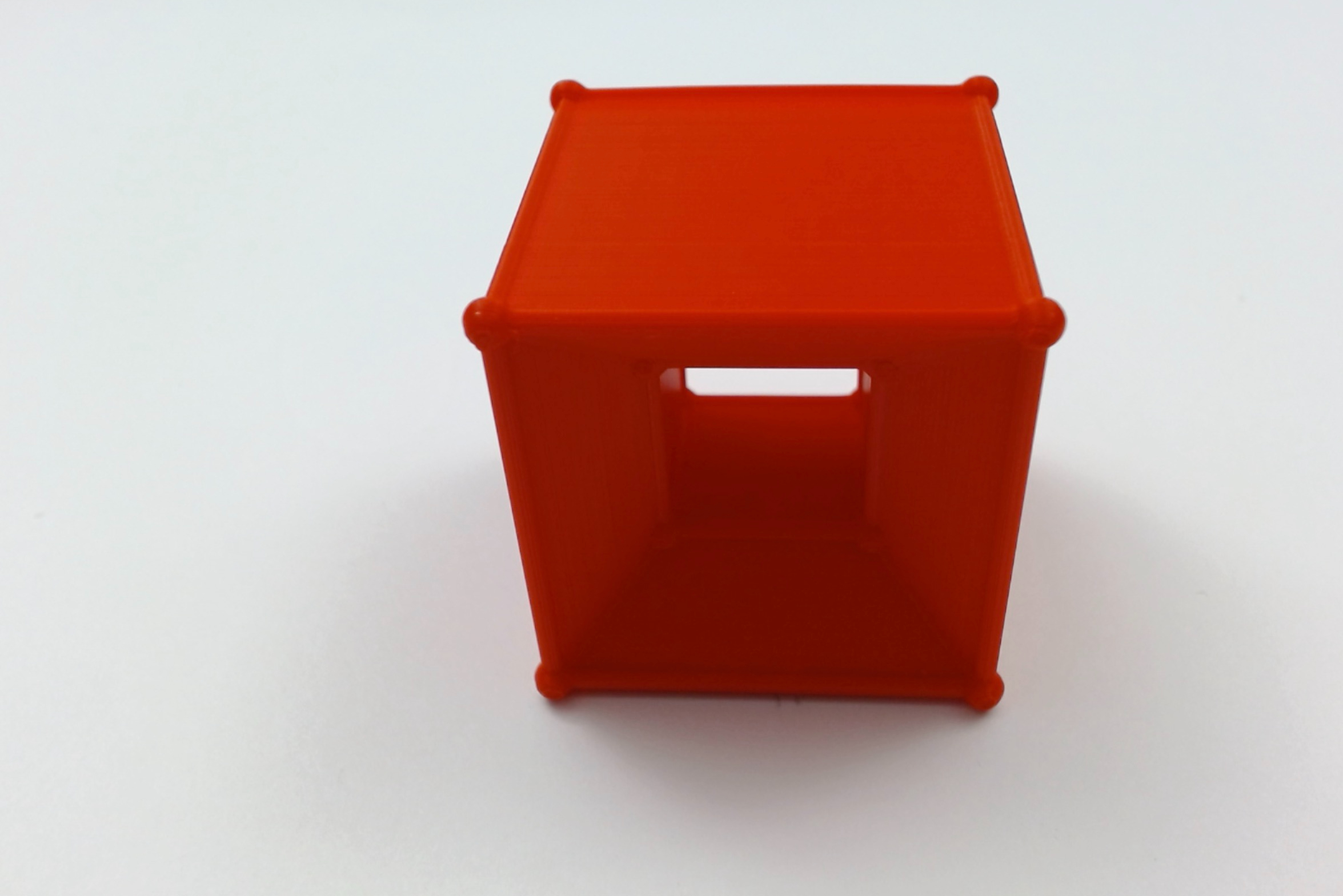}
    \quad
        \includegraphics[height=1.2in]{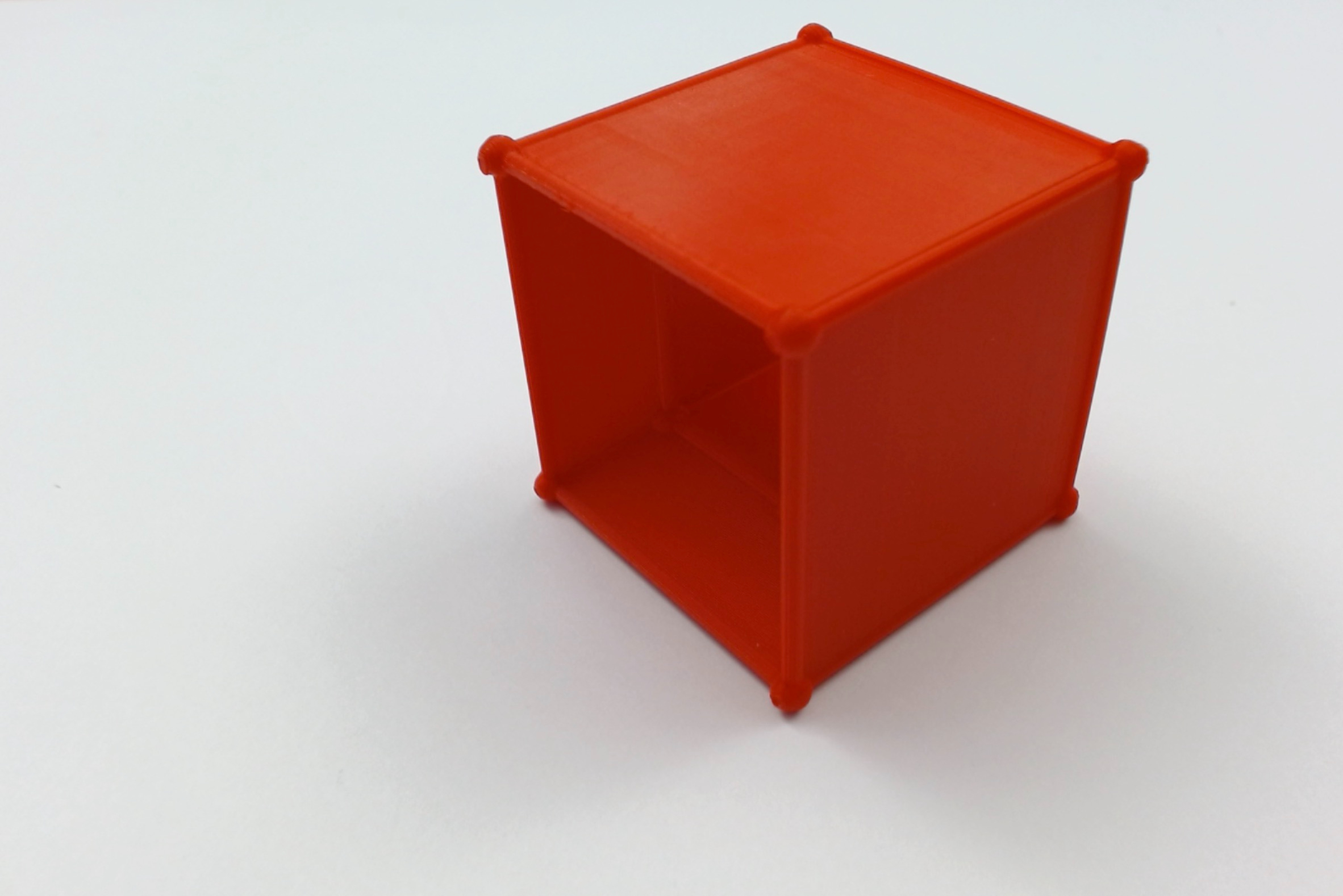}
\caption{Embedding of $\mathbb{T}^2$ in the tesseract.
The surface has 3 possible representatives.}
\label{print_T2}
\end{figure}

\begin{figure}
\centering
        \includegraphics[height=1.2in]{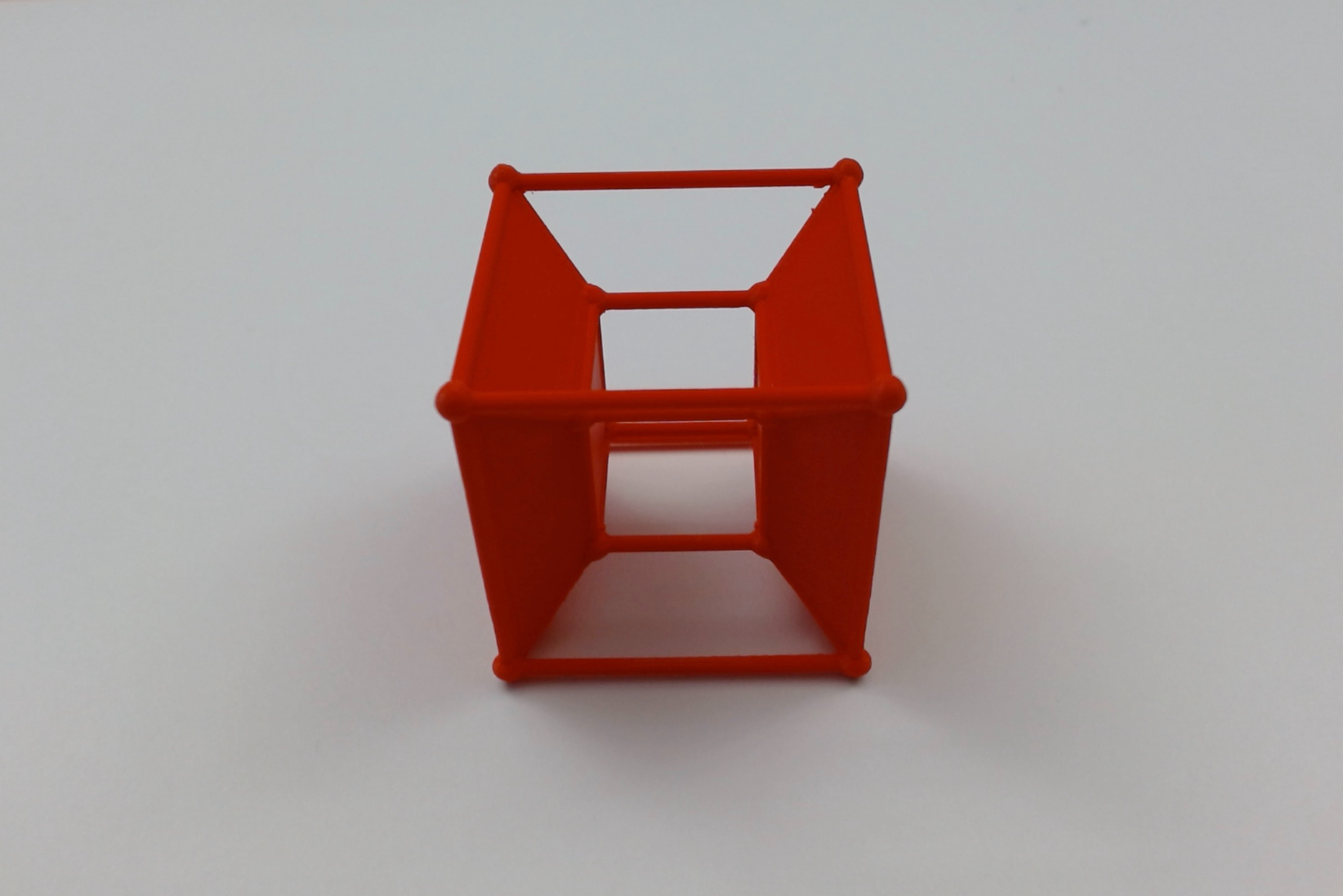}
    \quad
        \includegraphics[height=1.2in]{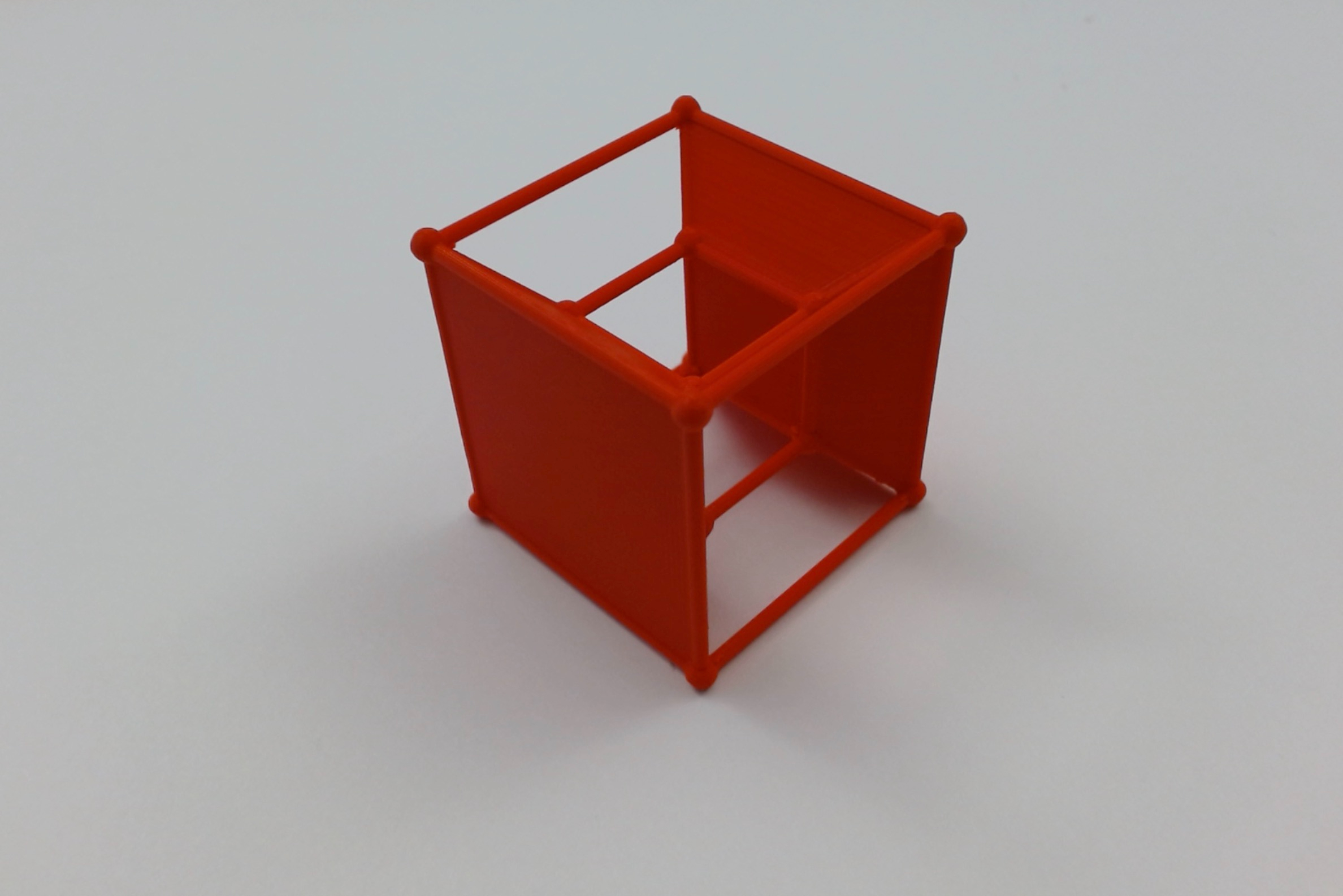}
\caption{Embedding of $\mathbb{S}^2 \coprod \mathbb{S}^2$ in the tesseract.
The surface has 4 possible representatives.}
\label{print_S2_S2}
\end{figure}

A two-dimensional cubical complex $\mathcal{C}$ is a \textit{cubical surface (without boundary)} if every point has a neighborhood $D$ that is homeomorphic to a two-dimensional disk. For the tesseract, an equivalent condition is that all edges of $\mathcal{C}$ are contained in exactly two faces. (More generally, there will be an extra condition at the vertices.) Using this characterization of a cubical surface on the tesseract, Aveni, Govc, and Rold\'an~\cite{AveniGovcRoldan} computed an exhaustive search of all cubical surfaces in the tesseract, finding 127 surfaces. Under the action of the symmetry group $B_4$, these fall into eight equivalence classes.
For each of these isomorphic types we have chosen a suitable representative for 3D printing. Figures~\ref{prints},~\ref{print_T2}, and~\ref{print_S2_S2} show two photographs of each.

Classifying and counting the possible configurations of some system is a very common problem in combinatorics. One of the most famous artist working in this direction was Sol LeWitt. Perhaps his closest work to ours is \emph{Variations of Incomplete Open Cubes}. Here LeWitt finds 122 connected subsets of the one-skeleton, up to orientation preserving symmetry (but leaving out the entire cube, and also subsets that lie in a single face of the cube).
We go up in dimension, selecting from the 24 faces of the tesseract rather than the 12 edges of the cube. 
This would result in far more configurations, but the extra constraint of the faces forming a closed surface greatly restricts the possibilities.

\section{A cubical surface with boundary}

A natural next question is to classify the \emph{cubical surfaces with boundary}.
Such a surface is a two-dimensional cubical complex with the property that any point has an open neighborhood homeomorphic to either a disk or a half-disk. 
This is a much weaker constraint on the possibilities than requiring that the surface be closed and we have yet to perform an exhaustive search for such surfaces. 
However, we found a cubical surface homeomorphic to a M\"obius strip. 
Note that the tesseract is the lowest dimensional cubical complex in which this is possible.
One possible realization of the M\"obius strip is given by the faces:
$$
\mathcal{C}:= \{ (*,*,1,1), (*,0,1,*), (0,0,*,*), (*,1,1,*), (*,1,*,0), (0,*,*,0) \}.
$$

\begin{wrapfigure}[11]{l}{0.30\textwidth}
\centering
\includegraphics[width=0.25\textwidth]{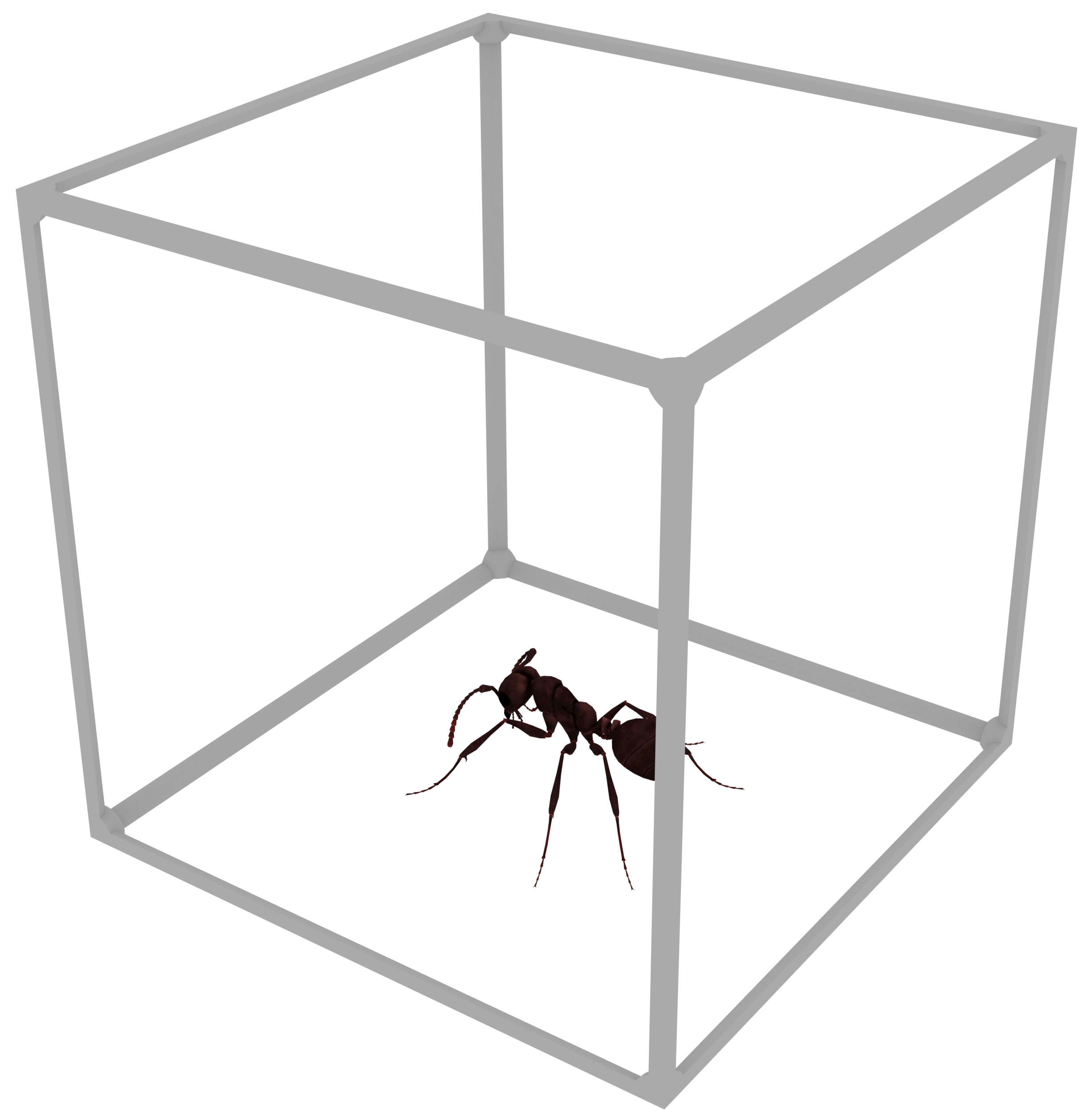}
\caption{Ant inside the cube.}
\label{Fig3dant}
\end{wrapfigure}

We could illustrate this M\"obius strip in the same style as our other surfaces, but we were inspired to create an homage to Escher's print \emph{M\"obius Strip II}.
To do this, we follow a similar procedure to the one implemented in~\cite{SegermanHart}. That is, we place a three-dimensional model of our ant inside of a cube in such a way that it is nicely standing on one of the faces of the cube (see Figure~\ref{Fig3dant}). Then we embed this cube as a cell of the tesseract. Finally, we project the faces with the ants back to three-dimensional space. In star notation, the eight three-dimensional cells of the tesseract are represented as 
$$
(0,*,*,*), (1,*,*,*), (*,0,*,*), (*,1,*,*), (*,*,0,*), (*,*,1,*), (*,*,*,0), (*,*,*,1).
$$
We may map our ant from the three-dimensional cube into each one of these cubes via
$$
(0,x,y,z), (1,x,y,z), (x,0,y,z), (x,1,y,z), (x,y,0,z), (x,y,1,z), (x,y,z,0), (x,y,z,1),
$$
selecting the permutation of the entries $x$, $y$, and $z$ in such a way that the ants are facing in the correct directions as they walk around the M\"obius strip.

\newpage

\begin{figure}[htb]
\centering
        \includegraphics[width=0.42\textwidth]{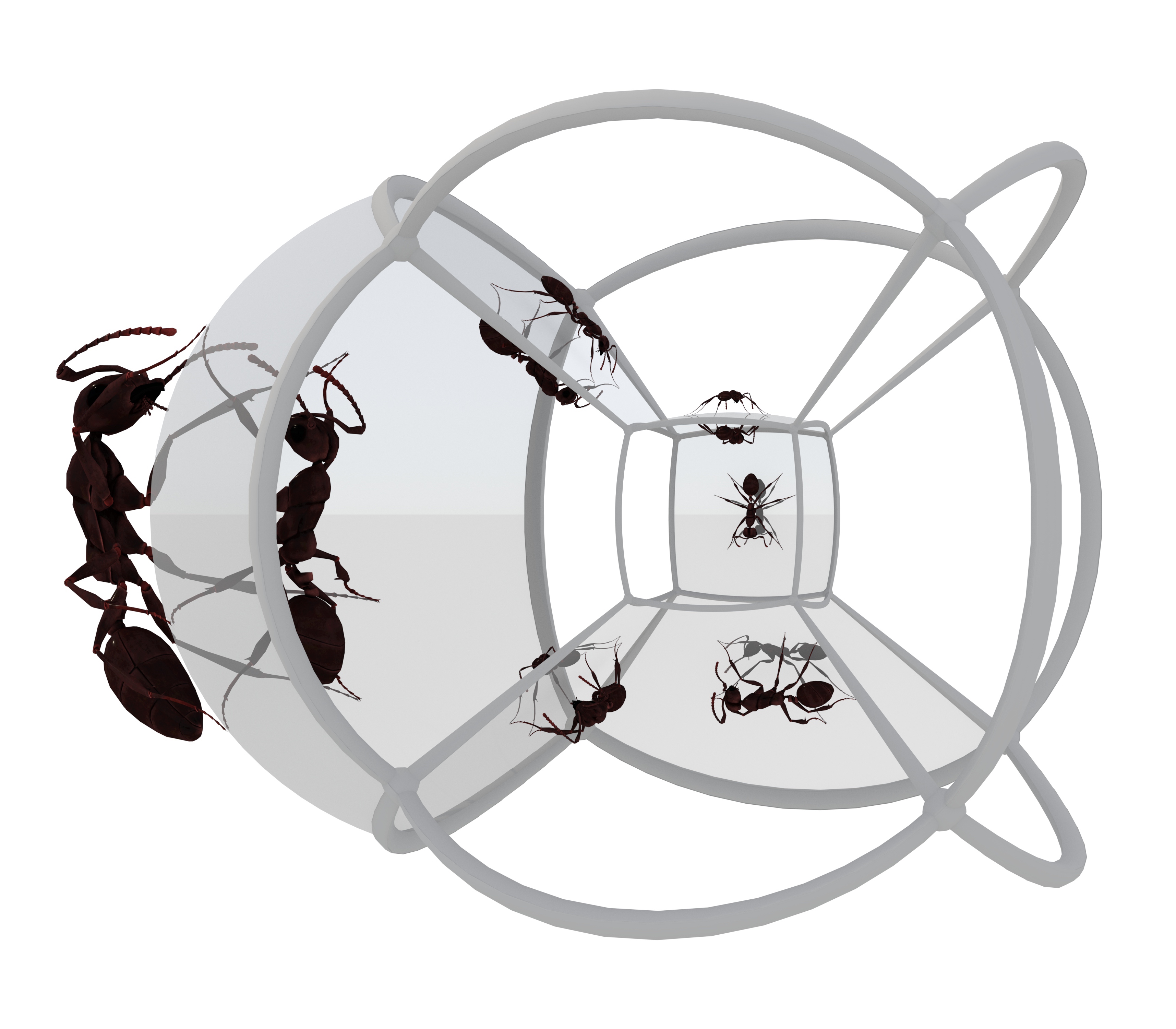}
\quad
        \includegraphics[width=0.48\textwidth]{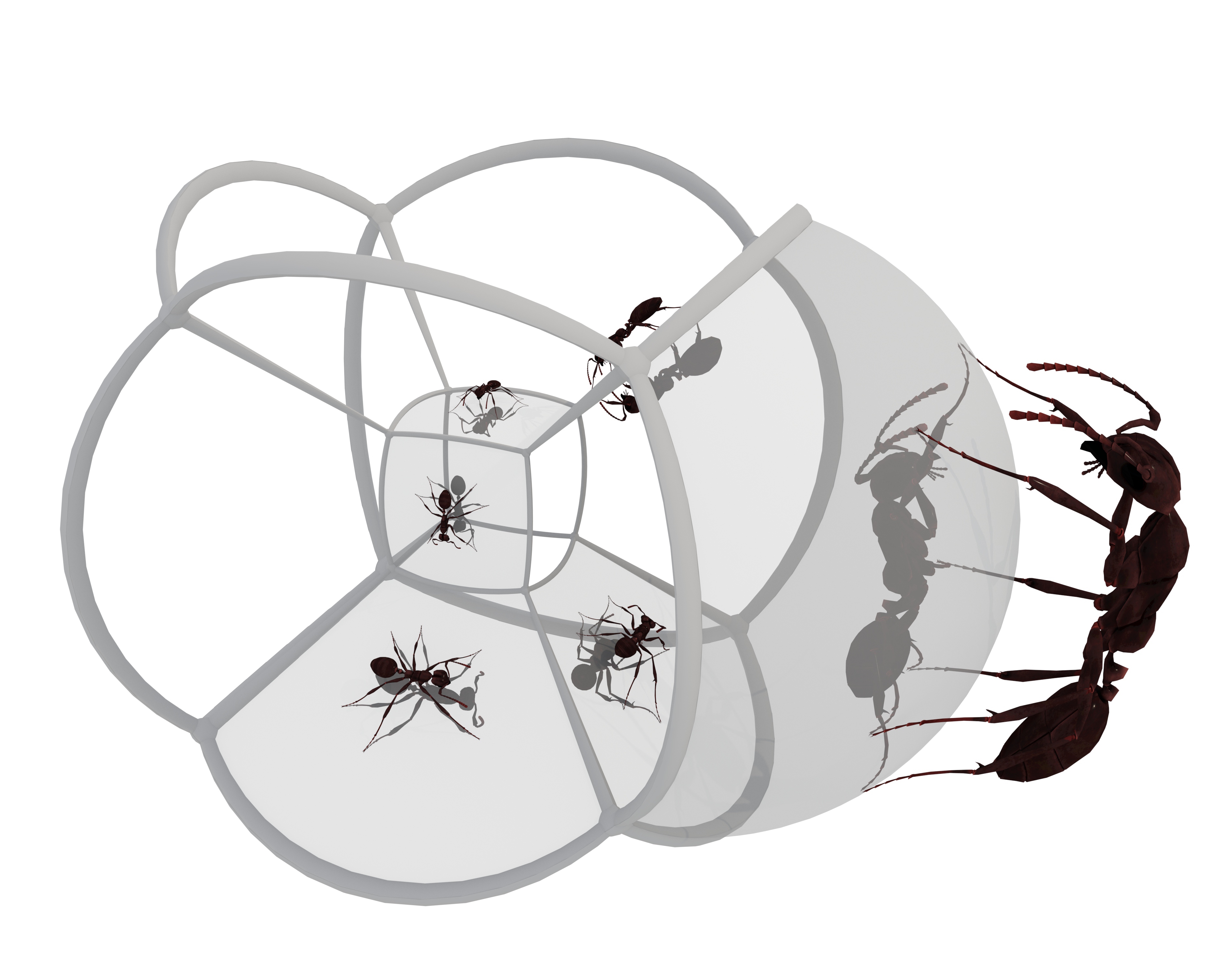}
\caption{Hyper-ants walking along the M\"obius strip.\\
Sketchfab link: \url{https://skfb.ly/oEMDt}}
\label{hypercubeants}
\end{figure}

To bring our ants back to three-dimensional space we use a composition of radial projection followed by stereographic projection to map the ant from the cube to the ``shadow'' of the tesseract. To achieve this, consider the unit vector $(0,0,0,1) \in \mathbb{R}^4$, this will be the ``north pole" of our hypersphere $S^3$, from which we will apply the stereographic projection.  To keep notation as clear as possible, consider the translation $\tau: \mathbb{R}^4 \rightarrow \mathbb{R}^4$ defined by $(x,y,z,w) \mapsto \left(x-\frac{1}{2},y-\frac{1}{2},z-\frac{1}{2},w-\frac{1}{2}\right)$, which will take the center of  $Q^4$ to the origin in $\mathbb{R}^4$. For any vertex $v \in Q^4_0$, the point $\tau(v)$ will have coordinates $\{\pm 1/2\}^4$. We use the radial projection $\rho :\mathbb{R}^4 \setminus \{(0,0,0,0)\} \rightarrow S^3$ defined by $v \mapsto \frac{v}{||v||}$. This maps any non-zero vector to its projection on the hypersphere. The stereographic projection mapping $\pi_p :\mathbb{R}^4 \setminus \{(0,0,0,1)\} \rightarrow \mathbb{R}^3$ defined by $(x,y,z,w) \mapsto \left(\frac{x}{1-w}, \frac{y}{1-w}, \frac{z}{1-w}\right)$
takes us back to $\mathbb{R}^3$.
Here we project from a point at the center of a cube of the tesseract. 
This means that the edges and faces of the tesseract are as small in the projection as possible (certainly we want to avoid any of these going through the projection point). Our ant crawls on, and is close to one of the faces of the cube. This means that it does not get too close to the center of its cube. Thus, even the ant in the ``outer'' cube of the projection is not too big.
Figure \ref{hypercubeants} shows the result, applied to the set of faces $\mathcal{C}$.


\section{Acknowledgements}
This project is supported by the Alexander von Humboldt Foundation. The third author was supported in part by National Science Foundation grant DMS-2203993. Thanks to [kyan0s, \url{https://sketchfab.com/kyan0s}] for the 3D model of the ant.

{\setlength{\baselineskip}{13pt} 
\raggedright				

} 
\end{document}